\documentclass{article}[11pt]
\usepackage{geometry}
\usepackage{enumerate}
\usepackage{latexsym,booktabs}
\usepackage{amsmath,amssymb, bm, tikz}
\usepackage{graphicx}
\usepackage{amsthm}
\usepackage{thmtools}
\usepackage{float}        
\usepackage{longtable}
\usepackage{multirow}
\usepackage{diagbox}

\usepackage{natbib}        
\usepackage{url}
\usepackage{textcomp}      
\usepackage{color}
\usepackage{subfiles}
\usepackage{mathtools}
\usepackage{multirow}

\usepackage{subcaption}

\usepackage{accents}

\usepackage[skip=2pt,font=footnotesize]{caption} 

\usepackage{siunitx} %
\usepackage[flushleft, online]{threeparttablex} %

\usepackage{tikz}
\usetikzlibrary{arrows.meta, positioning, decorations.markings, shapes}
\usepackage{pgfplots}
\pgfplotsset{compat=1.17}

\usepackage{parskip}

\usepackage{algorithmicx}
\usepackage[ruled]{algorithm}
\usepackage{algpseudocode}
\algdef{SE}[DOWHILE]{Do}{doWhile}{\algorithmicdo}[1]{\algorithmicwhile\ #1}%

\usepackage{cases}

\newcommand\hlight[1]{\tikz[overlay, remember picture,baseline=-\the\dimexpr\fontdimen22\textfont2\relax]\node[rectangle,fill=blue!50,rounded corners,fill opacity = 0.2,draw,thick,text opacity =1] {$#1$};}

\renewcommand{\arraystretch}{1.2}

\theoremstyle{definition}

\renewcommand\thmcontinues[1]{\textbf{Continued}}

\title{A graph-based heuristic for the non-stationary stochastic lot-sizing problem under penalty costs}

\usepackage{authblk}

\author[a]{Xiyuan Ma\thanks{Corresponding author: Xiyuan.Ma@ed.ac.uk}}
\author[a]{Roberto Rossi}
\author[a]{Thomas Welsh Archibald}
\affil[a]{Business School, University of Edinburgh, Edinburgh, United Kingdom}

\date{}

\begin{document}
\maketitle

\begin{abstract}

We present a relaxation‐and‐augmentation algorithm to compute near-optimal $(R,S)$ policy parameters for the single‐item single-stocking location non-stationary stochastic lot‐sizing problem under penalty costs. Our approach starts from a filtered state-space relaxation of the problem and then systematically repair infeasibilities through a repetitive state‐space augmentation procedure, until an optimal solution is obtained. We benchmark our approach against the state-of-the-art MILP-based cut generation approach for the problem; an extensive computational study based on a test bed from the literature demonstrates that our approach is computationally superior on medium and long horizon cases, and hence more scalable.  

\end{abstract}

\section{Introduction}
Effective inventory control is central to maintaining the efficiency of the supply chain in meeting fluctuating customer demands while controlling operational costs.
\citet{wemmerlov1989behavior} emphasise that inventory models must account for stochastic and dynamic environments to capture practical relevance and to further motivate research into non-stationary stochastic demand.

A key element of inventory control is lot sizing, whose aim is to determine \textit{replenishment timings} and associated \textit{order quantities} \citep{silver1981operations}.
Unfortunately, when demand is stochastic, replenishment timings and order quantities must be constantly revised in response to demand fluctuations; this leads to so-called setup-oriented and quantity-oriented system nervousness \citep{de1997nervousness,Tunc2013}. To address the longstanding issue of system nervousness in inventory control, \citet{bookbinder1988strategies} proposed the ``static-dynamic uncertainty'' strategy. In this strategy, one of these two key decisions, replenishment timings or order quantities, is made in a here-and-now fashion at the beginning of the planning horizon, while the other is made in a wait-and-see fashion. Depending on which decision is determined at the outset, the static-dynamic strategy results in two inventory policies: the ($R,S$) policy, in which replenishment timings\footnote{$R$ refers to the replenishment cycle, defined as the length of time that a single replenishment is intended to cover.} and associated order-up-to levels are fixed at the beginning of the planning horizon, and order quantities are determined, at each replenishment period, to raise the inventory level to the order-up-to level; and the ($s,Q$) policy, in which order quantities are fixed at the beginning of the planning horizon, while the timings of replenishments are determined dynamically, by means of the reorder threshold $s$.

In cost-based inventory optimisation, a fundamental trade-off exists between ordering frequently (higher ordering costs) and ordering infrequently in larger quantities (higher holding costs). At the core of this trade-off is the need to address unmet demand, which directly contributes to overall system cost.
There are two main modelling approaches aimed at controlling unmet demand.
The first is through service-level constraints \citep{Silver_Pyke_Peterson_2006}, which define performance targets such as a non-stockout probability ($\alpha$ service level), or the average fraction of demand that is fulfilled directly from available inventory ($\beta$ service level).
An alternative approach is to penalise unmet demand by means of a per-unit penalty cost. Finally, unmet demand can be treated as lost sales, or as backorders: in the former case,  demand is permanently lost; in the latter, all (complete backorders) or part of (partial backorders) demand is delayed and fulfilled as soon as a suitable replenishment arrives. This work focuses on the complete backorder case under penalty cost.

The aim of this work is to develop an efficient algorithm for computing near-optimal non-stationary ($R,S$) policy parameters\footnote{In a non-stationary ($R,S$) policy the lengths and order-up-to levels of replenishment cycles may vary with time.}. This problem has been widely studied in the literature, with a substantial body of work focusing on service level-based formulations. These models benefit from convenient analytical properties that ease the design of efficient solution methods. In particular, the optimal order-up-to level $S$ in each replenishment cycle can often be computed by inverting the cumulative distribution function (CDF) of demand \citep{tarim2004stochastic,tarim2006modelling,rossi2015piecewise}.
In contrast, the penalty cost formulation lacks such closed-form expressions. Determining $S$ in this setting requires solving a multi-period newsvendor problem for each candidate replenishment cycle, involving numerical evaluation of expected holding and backorder costs over multiple periods. Since the number of potential replenishment cycles grows quadratically with the planning horizon, the computational burden quickly becomes substantial, especially in non-stationary settings, in which demand varies over time.

To address these challenges, we propose a novel, computationally efficient, graph-based heuristic that formulates the non-stationary $(R,S)$ penalty cost model as a shortest-path problem with state-space filtering.
Rather than constructing the full graph in advance, we apply dominance-based pruning during graph generation, leveraging triangle inequality violations to discard unpromising paths. This dramatically reduces the number of arcs that must be evaluated and enables fast computation even for long planning horizons.

We make the following contributions to the literature.
\begin{itemize}\setlength\itemsep{-0.5em}
  \item We develop a relaxation‐and‐augmentation framework to compute near-optimal ($R,S$) policy parameters for the single-item single-stocking location non-stationary stochastic lot sizing problem under penalty cost. We start from a filtered state-space relaxation of the problem and then systematically repair infeasibilities through a repetitive state‐space augmentation procedure, until an optimal solution is obtained. 
  \item We conduct an extensive computational study on 1620 benchmark instances to compare our approach against the state-of-the-art cut generation method of \citet{tunc2018extended}. Our study demonstrates that our approach is computationally superior on medium and long horizon cases.
\end{itemize}

The remainder of the paper is organised as follows. 
Section~\ref{sec:LR} surveys related literature. 
Section~\ref{sec:RS-SP} provides the problem statement and the associated stochastic programming formulation.
Section~\ref{sec:DE} introduces a deterministic non-linear mathematical programming formulation for the problem. Section~\ref{sec:SPF} discusses a shortest-path relaxation for the problem. Section~\ref{sec:RSalgorithm} presents our state space augmentation algorithm.
Section~\ref{sec:computation} reports on our computational study, in which we compare our approach against the state-of-the-art cut-generation method of \cite{tunc2018extended}. 
Section~\ref{sec:conclusion} concludes the paper and outlines directions for future research.

The notation used throughout the paper is summarised in Appendix~\ref{sec:App.notations}.

\section{Literature review}\label{sec:LR}
We survey the literature that investigates scalable and effective strategies for computing non-stationary ($R,S$) policy parameters.

A considerable body of research has addressed the determination of non-stationary ($R, S$) policy parameters under service-level constraints. \citet{tarim2004stochastic} introduced a mixed-integer programming (MIP) formulation for the non-stationary ($R, S$) policy under $\alpha$ service level constraints. This formulation was later extended to alternative service level definitions, including the $\alpha_c$ by \citet{tempelmeier2007stochastic}, and the $\beta$ and $\beta^{cyc}$ levels by \citet{rossi2015piecewise}. \citet{tunc2014reformulation} decouple replenishment cycles using a network-flow-based MIP structure, resulting in tighter linear relaxations.

In parallel to service-level-based models, \citet{tarim2006modelling} proposed a penalty cost-oriented MIP formulation using piecewise linear approximations of expected shortage costs. The approximation method, later refined by \citet{rossi2014piecewise}, allows for adjustable accuracy and provides piecewise linearisation parameters tailored to the normal distribution. To enhance computational performance, \citet{tunc2018extended} generalised this approach with a dynamic cut generation algorithm, incorporating the piecewise cost structure from \citet{rossi2015piecewise}.

A common assumption in these models is that \textit{negative orders are disallowed}, in other words, excess stock cannot be returned if realised demand is significantly lower than expected. In practice, this implies that inventory carried over exceeds the order-up-to level, but cannot be “sold back” to the supplier. \citet{tarim2011efficient} relaxed this constraint and proposed, under service level constraints, a computationally efficient MILP-based heuristic that solves a relaxed shortest-path formulation. Even when infeasible, the relaxed solution provides tight lower and upper bounds on the optimal cost; this relaxation is then integrated into a branch-and-bound framework. In parallel, building upon the framework in \citep{Christofides1981}, \citet{rossi2011state} developed a state-space augmentation method that repairs the relaxed state-space-graph in order to produce a complete and compact representation of the original problem under service level constraints. \citet{ozen2012static} analysed the optimality of base-stock policies under both service-level constraints and penalty cost formulations, and proposed two heuristics: the approximation and the relaxation heuristic; in the latter heuristic, non-negativity constraints on order sizes are relaxed.

The majority of existing studies on the non-stationary $(R, S)$ policy focus on service-level models, while relatively few have addressed the penalty cost formulation. Existing studies focusing on the backorder penalty cost setting typically present MIP-based heuristics, while no specialised efficient algorithms akin to \citet{rossi2011state} have been presented so far. This study addresses this gap by introducing a novel, graph-based heuristic algorithm to efficiently compute near-optimal non-stationary $(R, S)$ policy parameters under penalty costs.

\section{Problem statement and Stochastic Programming formulation}\label{sec:RS-SP}
We consider the single-item single-stocking location non-stationary stochastic lot-sizing problem over a planning horizon comprising $T$ periods.
Demand $d_t$ in period $t\in T$ is a random variable with cumulative distribution function $\Phi_t(\cdot)$; demands in different periods are independently distributed.
We consider both a fixed ordering/production setup cost component $K$, and a proportional ordering/production cost component $z$.\footnote{Note that in a periodic review inventory/production setting the fixed ordering/production setup cost associated with a replenishment is incurred {\em independent of the batch size} \citep[][p. 44]{axsater2015inventory}; this is due to the fact that $K$ represents expenses such as setup and learning costs, administrative tasks, transportation and material handling, and equipment downtime associated with the  production/inventory system.} Negative orders are not allowed, that is $Q_t\geq 0$.
At the end of each period, a  holding cost $h$ is incurred for each unit of inventory carried to the next period or beyond the end of the planning horizon. Without loss of generality, we assume replenishments are delivered instantaneously.
We operate under complete backorders: unmet demand is satisfied as soon as a suitable replenishment arrives.
In each period, a penalty cost $b$ is incurred for any unit of unmet demand.
The objective is to minimise the expected total cost, which comprises fixed/variable ordering costs, holding costs, and penalty costs.

Let $R(i,j)$ denote a replenishment cycle in which an order/production setup takes place in period $i$ to cover demand until period $j$, and let $S_i$ denote the corresponding order-up-to level. Recall that, under an ($R,S$) policy, replenishment cycles and corresponding order-up-to levels are determined simultaneously at the beginning of the planning horizon, while actual order quantities are decided at the beginning of each replenishment cycle to raise inventory to the corresponding order-up-to level \citep{bookbinder1988strategies}.

Accordingly, we formulate the problem of computing non-stationary ($R,S$) policy parameters as a two-stage stochastic programming model.
Let $\delta_t$ be a binary decision variable that equals 1 if an order/setup occurs at the beginning of period $t$, and 0 otherwise.
These decision variables represent the replenishment timings so that, if the replenishment cycle $R(i,j)$ is included in an optimal policy, $\delta_i = 1$ and $\delta_k = 0$ for $1 < k \le j$.
%
Let $I_t$ be a random variable that represents the closing inventory level of period $t$, where $I_0$ denotes the initial inventory level.
Finally, let $Q_t$ be a non-negative random variable representing the order quantity in period $t$.

The model ($P$) is as follows,


\begin{alignat}{2}
C_1(I_0) &= \min_{\substack{\{\delta_t\}_{t=1}^T \\\\ \{S_t\}_{t=1}^T}} \quad
\mathrm{E}\Big[\displaystyle\sum_{t=1}^T
\Big(\delta_t K + zQ_t + h\max(I_t, 0) + b\max(-I_t,0)\Big)\Big] \label{eq:SPobj}\\
\mbox{s.t.}\quad
& \text{for }t=1,\ldots, T,\notag\\
&I_t = I_0+\sum_{i=1}^t(Q_i-d_i),\label{eq:SPconservation}\\
&\delta_t = 0 \rightarrow Q_t = 0, \label{eq:SP-Q-1}\\
&\delta_t = 1 \rightarrow Q_t = \max(S_t - I_{t-1},0), \label{eq:SP-Q-2}\\
&I_t\in\mathbb{R}, \quad \delta_t\in\{0,1\}, \quad Q_t\geq 0, \quad S_t\geq 0, \label{eq:SPdomain}
\end{alignat}

where $\rightarrow$ denotes the implication constraint \citep[see e.g.][]{cplex_user_manual} and $\mbox{E}$ denotes the expected value operator.
In this model, first stage decision variables are the order/production setup timings $\delta_t$ and associated order-up-to levels $S_t$; second stage recourse variables are the order quantities $Q_t$ and the closing inventory levels $I_t$. From this decision structure, it is clear that the (fixed) cost entailed by setting any $\delta_t=1$ is incurred at the beginning of the planning horizon, before any demand is realised, and before any realised order/production quantity is known; in other words, we are incurring fixed costs to {\em schedule} a future order/production setup; this reflects the nature of inventory reviews in the ($R,S$) inventory policy. The objective \eqref{eq:SPobj} is to minimise the expected total cost; constraint \eqref{eq:SPconservation} enforces inventory conservation constraints; \eqref{eq:SP-Q-1} and \eqref{eq:SP-Q-2} define $Q_t$ as the order quantity needed to raise inventory to the order-up-to level $S_t$, if an order is scheduled in period $t$, and zero otherwise; note that if, at a review period, inventory $I_{t-1}$ is equal or exceeds $S_t$, then the order quantity is zero.

\section{A deterministic non-linear reformulation}\label{sec:DE}
The two-stage stochastic programming model presented in the previous section can be reformulated, under mild assumptions \citep[see e.g.][]{tarim2006modelling,tunc2018extended}, by a deterministic non-linear model $(D)$ as follows,

\begin{alignat}{2}
C_1(I_0) &= \min_{\substack{\{\delta_t\}_{t=1}^T \\\\ \{S_t\}_{t=1}^T}} \quad
\displaystyle\sum_{t=1}^T
\Big(\delta_t K + h H_t + b B_t\Big) + z \tilde{I}_t + \underbrace{z \sum_{t=1}^T \tilde{d}_t - z I_0}_{\text{constant}} \label{eq:DEobj}\\
\mbox{s.t.}\quad & \mbox{for } t=1,2,\ldots, T\notag\\
&\delta_t = 0 \rightarrow \tilde{I}_t + \tilde{d}_t - \tilde{I}_{t-1} = 0, \label{eq:DE-Q-1}\\
&\delta_t = 1 \rightarrow S_t = \tilde{I}_t + \tilde{d}_t, \label{eq:DE-Q-3}\\
&\tilde{I}_t + \tilde{d}_t - \tilde{I}_{t-1} \geq 0, \label{eq:DE-Q-2}\\
&\sum_{j=1}^tP_{jt}=1\label{eq:DE-Pjt-1}\\
& \tilde{I}_t = H_t - B_t \label{eq:DE-Ht-Bt} \\
&P_{jt}=1 \rightarrow H_t = \mbox{E}[\max(S_j-d_{jt},0)] \label{eq:DE-Ht} & j=1,2,\ldots,t\\
&P_{jt}=1 \rightarrow B_t = \mbox{E}[\max(d_{jt}-S_j,0)]\label{eq:DE-Bt} & j=1,2,\ldots,t\\
&P_{jt}\in\{0,1\}, &j=1,2,\ldots,t\\
&\delta_t\in\{0,1\},\\
&S_t\geq 0, \quad \tilde{I}_t\in\mathbb{R}, \quad H_t,B_t\geq 0.
\end{alignat}

where $d_{jt}=d_j+\ldots+d_t$, $\tilde{d}_t$ is the expected demand in period $t$, $\tilde{I}_t$ denotes the expected inventory level, $B_t$ represents the expected negative inventory and $H_t$ the expected positive inventory at the end of period $t$.
To obtain Eq.\eqref{eq:DEobj} from Eq.\eqref{eq:SPobj}, we push the expectation inside the summation --- since $\delta_t$ is a here-and-now decision fixed at the beginning of the planning horizon, it follows that $\mathbb{E}[\delta_t K] = \delta_t K$. Similarly, since $S_t$ is also fixed at the outset under the $(R,S)$ policy, the expected holding and backorder cost terms $hH_t$ and $bB_t$ are obtained via constraints \eqref{eq:DE-Ht} and \eqref{eq:DE-Bt}, after noting that $I_t = S_j - d_{jt}$ for all $t$. Finally, the proportional ordering cost term follows from observing that $\sum_{t=1}^T Q_t = \sum_{t=1}^T \tilde{d}_t + \tilde{I}_T - I_0$.
Note that $\mbox{E}[\max(d_{jt}-S_j,0)]$ is the first order loss function; $P_{jt}$ is a binary variable that takes value 1 if the most recent replenishment before period $t$ occurred in period $j$. Finally, note that in Eq.\eqref{eq:DE-Q-2} the term $\tilde{I}_t + \tilde{d}_t - \tilde{I}_{t-1} $ represents the expected order quantity in period $t$; and if $\delta_t=1$, it follows that $S_t=\tilde{d}_t+\tilde{I}_t$. The loss function in this model can be piecewise linearised, and the model can be solved using off-the-shelf mathematical programming solvers \citep{rossi2015piecewise}. The optimality gap of these approximations has been thoroughly investigated in the literature \citep[see e.g.][]{tunc2018extended, DuralSelcuk2020}, therefore, in {\em this work we will focus on developing a computationally efficient solution method to produce optimal solutions to model $D$}.

\paragraph{Numerical example.} We consider the following instance: demands in different periods are independently distributed and follow the normal distribution with means $\tilde{d}_t = \{100,125,25,40,30\}$ and coefficient of variation $\sigma_t/\tilde{d}_t = 0.3$; other problem parameters are $b=19$, $h=1$, $K=60$ and $z=0$. By solving the instance using the approach of \citep{rossi2015piecewise} with 11 segments, we obtain the replenishment plan illustrated in Figure \ref{fig:replenishment_plan_rossi2015}, and estimated lower (478) and upper (494) bounds for its expected total cost. Note that replenishment periods and associated order-up-to levels are identified by red dots.
\begin{figure}
\centering
\resizebox{.43\textwidth}{!}{
\begin{tikzpicture}
    \begin{axis}[
        xlabel={Period},
        ylabel={inventory level},
        xmin=0, xmax=5,
        ymin=-50, ymax=220,
        xtick={0,1,...,5},
        xticklabels=\empty,
        extra x ticks={0.5,1.5,2.5,3.5,4.5,5.5},
        extra x tick labels={1,2,3,4,5},
        extra x tick style={major tick length=0pt},
        ytick={0,50,...,200},
        legend pos=north east,
        ymajorgrids=true,
        grid style=dashed
    ]
    \addplot[
        color=blue,
        mark=x,
        ]
        coordinates {
        (0,142)(1,142-100)(1,203)(2,203-125)(3,203-125-25)(3,91)(4,91-40)(5,91-40-30)
        };
    \addplot[
        color=red,
        only marks,
        mark=*,
        ]
        coordinates {
        (0,142) (1,203) (3,91)
        };
    \end{axis}
\end{tikzpicture}
}
\caption{Near-optimal replenishment plan for our running example obtained via the MILP approach in \citep{rossi2015piecewise}}
\label{fig:replenishment_plan_rossi2015}
\end{figure}

\section{A shortest path relaxation}\label{sec:SPF}
By following \citep{tarim2011efficient, ozen2012static},  we relax expected order quantity nonnegativity constraints in model $D$, and obtain a relaxed model $R$, which we reformulate as a shortest-path problem.

We relax constraint \ref{eq:DE-Q-2} by allowing negative expected order quantities, i.e. $\tilde{I}_t + \tilde{d}_t - \tilde{I}_{t-1}<0$, which represents inventory adjustments in which excess stock is shipped back to the warehouse. This relaxation eliminates dependencies among replenishment cycles. In practice, this means that, given a cycle $R(i,j-1)$ and its order-up-to level $S_i$, its expected total cost can be computed, independently from other cycle order-up-to levels and expected total costs, by considering the following multiperiod newsvendor problem
\[C_{(i,j)}(S)=K+\mbox{E}\left[\sum_{k=i}^{j-1} h\max(S-d_{ik}, 0) + b\max(d_{ik}-S,0)\right].\]
By following \cite[][Eq. 3]{Askin1981}, the optimal $S_i$ for a given cycle $R(i,j-1)$ can be computed as the value $S$ that satisfies condition
\[\sum_{k=i}^{j-1} \Phi_t(S) = (j-i)\frac{b}{b+h}.\]
As remarked by Askin, this equation is fairly easy to solve since many terms in the sum will be approximately one for $k<j$. The expected total cost of cycle $R(i,j-1)$ is then obtained by plugging $S_i$ into $C_{(i,j)}(\cdot)$.

We now describe the construction of the directed acyclic graph (DAG) that we will use to model our relaxed problem.
Let $G\triangleq\langle V,E\rangle$ be a DAG, in which
nodes $V={1,\dots,T,T+1}$ represent replenishment timings plus a dummy terminal node $T+1$;
and an arc $(i,j) \in E$, with $1\leq i<j\leq T+1$, represents replenishment cycle $R(i,j-1)$ with associated cost $C_{i,j}(S_i)$.

Consider the shortest path from node 1 to node T+1, and recall that each arc $(i,j)$ represents a replenishment cycle $R(i,j-1)$ with associated cost $C_{i,j}(S_i)$; by construction, since replenishment cycles are now independent, the cost $C^*$ of the shortest path is the cost of an optimal policy for problem $R$.


It is worth observing that $C_{(i,j)}$ does not appear to capture the proportional ordering cost $z$. To consider the proportional ordering cost, the graph must be slightly amended, so that the cost of any arc $(i,T+1)\in\mathcal{E}$ is equal to $C_{(i,T+1)}(S_i)+z\tilde{I}_T$ --- this accounts for the proportional unit cost associated with the expected amount of inventory that remains at the end of the planning horizon; finally, we must add the expected total proportional ordering cost over the planning horizon, which is constant and equal to $(\tilde{d}_1+\ldots +\tilde{d}_T)z$, to the cost $C^*$ of the shortest path after we have computed it (see Eq. \ref{eq:DEobj}). Similarly, it is possible to slightly amend the graph structure to account for a positive initial inventory \citep[see][]{ozen2012static}.

To summarise, the arcs belonging to the shortest path uniquely determine a set of replenishment timings, and in turn these timings partition the entire planning horizon into non-overlapping replenishment cycles $(i,j)$ whose order-up-to level $S_i$ can be determined independently, based on aggregate demand $d_{i,\ldots,j-1}$.
The presence of negative expected order quantities violates Eq.\eqref{eq:DE-Q-2} in problem $D$ and can be immediately checked once the shortest path is obtained. 
If no negative expected order quantity is observed, the replenishment plan that corresponds to the shortest path in the relaxed problem is optimal for the original problem. 
This optimality relies on the convexity of the multi-period newsvendor cost function $C_{ij}(\cdot)$, which means that each replenishment cycle cost can be minimised independently.
Conversely, if negative expected order quantities are observed, this replenishment plan is not feasible for problem $D$. To address this latter scenario, in the next section, we describe a state-space augmentation algorithm that iteratively repairs these infeasibilities.

\paragraph{Numerical example.}

We relax constraint \ref{eq:DE-Q-2} and allow negative expected order quantities for our running numerical example. We obtain the state-space shown in Figure \ref{fig:shortestPathExample}. The shortest path from node 1 to node $T+1$ has a cost of 475.4 and corresponds to the replenishment plan illustrated in Figure \ref{fig:relaxed_plan_example}, which features an expected negative order at the beginning of period 3.

\begin{figure}
\begin{minipage}[b]{0.55\columnwidth}
\begin{center}
\resizebox{1\textwidth}{!}{%
\begin{tikzpicture}[
  roundnode/.style={circle, draw, minimum size=0.6cm, inner sep=1pt},
  dashedarrow/.style={dashed, ->, -Latex},
  solidarrow/.style={->, -Latex},
  node distance=1.2cm and 1.5cm
]
  \node[roundnode, minimum size=1cm] (1) {1};
  \node[roundnode, right=of 1, minimum size=1cm] (2) {2};
  \node[roundnode, right=of 2, minimum size=1cm] (3) {3};
  \node[roundnode, right=of 3, minimum size=1cm] (4) {4};
  \node[roundnode, right=of 4, minimum size=1cm] (5) {5};
  \node[roundnode, right=of 5, minimum size=1cm] (6) {6};

  \draw[solidarrow, thick] (1) to [bend left=40] node[below] {\bf 121}  (2);
  \draw[solidarrow, thick] (2) to [bend left=40] node[below] {\bf 137} (3);
  \draw[solidarrow, thick] (3) to [bend left=40] node[below] {\bf 75.4} (4);
  \draw[solidarrow] (4) to [bend left=40] node[below] {84.7} (5);
  \draw[solidarrow] (5) to [bend left=40] node[below] {78.5} (6);

  \draw[solidarrow] (1) to [bend left=40] node[below] {480} (4);
  \draw[solidarrow] (2) to [bend left=40] node[below] {356} (5);
  \draw[solidarrow] (3) to [bend left=40] node[below] {238} (6);

    \draw[solidarrow] (1) to [bend right=40] node[below] {653} (5);
    \draw[solidarrow] (2) to [bend right=40] node[below] {500} (6);

  \draw[solidarrow] (1) to [bend right=40] node[below] {353} (3);
  \draw[solidarrow] (2) to [bend right=40] node[below] {225} (4);
  \draw[solidarrow] (3) to [bend right=40] node[below] {149} (5);
  \draw[solidarrow, thick] (4) to [bend right=40] node[below] {\bf 142} (6);

  \draw[solidarrow] (1) to [bend left=40] node[below] {837} (6);

\end{tikzpicture}
}
\caption{Graph structure representation of the relaxed problem for our numerical example, the shortest path is in bold, arcs are labelled with the respective cycle costs}
\label{fig:shortestPathExample}
\end{center}
\end{minipage}
\hfill
\begin{minipage}[b]{0.35\columnwidth}
\begin{center}
\centering
\resizebox{0.9\textwidth}{!}{
\begin{tikzpicture}
    \begin{axis}[
        xlabel={Period},
        ylabel={inventory level},
        xmin=0, xmax=5,
        ymin=-50, ymax=200,
        xtick={0,1,...,5},
        xticklabels=\empty,
        extra x ticks={0.5,1.5,2.5,3.5,4.5,5.5},
        extra x tick labels={1,2,3,4,5},
        extra x tick style={major tick length=0pt},
        ytick={0,50,...,200},
        legend pos=north east,
        ymajorgrids=true,
        grid style=dashed,
    ]
    \addplot[
        color=blue,
        mark=x,
        ]
        coordinates {
        (0,149)(1,149-100)(1,186)(2,186-125)(2,37)(3,37-25)(3,89)(4,89-40)(5,89-40-30)
        };
    \addplot[
        color=red,
        only marks,
        mark=*,
        ]
        coordinates {
        (0,149) (1,186) (2,37) (3, 89)
        };
    \end{axis}
\end{tikzpicture}
}
\caption{Replenishment plan for our running example obtained via the shortest path relaxation}
\label{fig:relaxed_plan_example}
\end{center}
\end{minipage}
\end{figure}

\section{A state space augmentation algorithm}\label{sec:RSalgorithm}
In presence of negative expected order quantities, we employ a  State Space Augmentation (SSA) algorithm to repair infeasibilities. This section details the methodology used to iteratively augment the state space and obtain an optimal solution to problem $D$.

\subsection{Generation and filtering of the state space graph}\label{sec:state_space_filtering}
When dealing with very large instances, it is important to ensure the generated state space graph is as compact as possible from the onset. To ensure this, we concurrently apply state-space filtering when the state space graph is initially generated, in order to minimise the size of the graph produced.

The shortest path algorithm is built upon the principle of optimality \citep{BellmanRichard1957Dp}, which states that any subpath of an optimal path is itself optimal. By leveraging this principle, we can prevent the generation of redundant arcs to reduce the problem size without affecting the optimal solution to the relaxed problem. However, while doing so, one must be careful to not prune arcs that may become necessary at a later stage, when infeasibilities are repaired.

Algorithm~\ref{alg:digraph} presents the pseudocode for the digraph generating procedure.
The algorithm proceeds as follows: the graph is initialized with $T+1$ nodes, where node $i$ corresponds to period $i$ (lines \ref{line:generation_start}--\ref{line:generation_end}). For every period $j$ (from $1$ to $T$), the algorithm loops over all earlier periods $i$ (with $i \geq i\_start$) and constructs a multiperiod newsvendor model $m$ for the subproblem spanning periods $i+1$ to $j$; then an edge from node $i$ to node $j$ is inserted into the graph with attributes: optimal order quantity $S$, expected cycle cost $C_{i,j}(S)$, and residual expected cycle closing inventory level $\tilde{I}$ (lines \ref{line:mpnb_start}--\ref{line:mpnb_end}). State-space filtering eliminates arcs that violate the triangle inequality: for each candidate intermediate node $k$, the algorithm checks whether the arc from $i\_start$ to $k$, when combined with an arc from $k$ to $j$, results in a lower cost than a direct arc from $i\_start$ to $j$. Additionally, it verifies that the expected closing inventory level in the first arc does not exceed the order-up-to level in the following transition. If these conditions hold and no previous negative orders were detected, $i\_start$ is incremented, effectively skipping arcs that cannot be part of an optimal solution (lines \ref{line:ssf_start}--\ref{line:ssf_end}). Finally, the directed acyclic graph encapsulating the state space and feasible transitions with their costs is returned. This digraph forms the basis for the subsequent formulation of the shortest path relaxation and state-space augmentation strategy.

\paragraph{Numerical example.} We apply the generation and filtering of the state space graph to our running numerical example; this leads to the reduced state space graph illustrated in Figure \ref{fig:shortestPathExample_filtered}. Observe that, in principle, arc (2,4) is dominated by the sequence $\langle (2,3), (3,4) \rangle$, as $225>137+75.4$; however, since an expected negative order is scheduled between cycles (2,3) and (3,4), arc (2,4) is retained in the filtered graph, as this arc may become candidate optimal when infeasibilities are repaired.
\begin{algorithm}[H]
\caption{Create Directed Graph for SSA}\label{alg:digraph}
\begin{algorithmic}[1]
\Procedure{create\_di\_graph}{}
    \State $T \gets $ number of periods \label{line:generation_start}
    \State $G \gets$ new directed graph
    \For{$i = 1$ to $T+1$}
        \State Add node $i$ to $G$ with label $i$ \label{line:generation_end}
    \EndFor
    \State $i\_start \gets 0$ \label{line:mpnb_start}
    \For{$j = 1$ to $T+1$}
        \For{$i = i\_start$ to $j-1$}
            \State Construct a multiperiod newsvendor model $m$ for $R(i,j-1)$
            \State $S \gets$ optimal order-up-to level for $R(i,j-1)$
            \State Add edge $(i, j)$ to $G$ with attributes $S$, $C=C_{i,j}(S)$, $\tilde{I}$, and $m$ \label{line:mpnb_end}
        \EndFor
            \For{$k = i\_start+1$ to $j-1$} \Comment{state-space filtering} \label{line:ssf_start}
                \State $no\_prev\_neg \gets $ \textbf{True}
                \For{$z = 0$ to $i\_start-1$}
                    \If{edge $(z,i\_start)$ exists \textbf{and} $(z,i\_start).\tilde{I} > (i\_start,k).S$}
                        \State $no\_prev\_neg \gets $ \textbf{False}
                        \State \textbf{break}
                    \EndIf
                \EndFor
                \If{$(i\_start,k).C + (k,j).C < (i\_start,j).C$} \Comment{triangle inequality}
                	   \If{$(i\_start,k).\tilde{I} \leq (k,j).S$ \textbf{and} $no\_prev\_neg$} \Comment{negative orders}
                       \State $i\_start \gets i\_start + 1$
                       \State \textbf{break} from the current filtering loop \label{line:ssf_end}
                    \EndIf
                \EndIf
            \EndFor
    \EndFor
    \State \Return $G$
\EndProcedure
\end{algorithmic}
\end{algorithm}

\begin{figure}
\centering
\resizebox{0.8\textwidth}{!}{%
\begin{tikzpicture}[
  roundnode/.style={circle, draw, minimum size=0.6cm, inner sep=1pt},
  dashedarrow/.style={dashed, ->, -Latex},
  solidarrow/.style={->, -Latex},
  node distance=1.2cm and 1.5cm
]
  \node[roundnode, minimum size=1cm] (1) {1};
  \node[roundnode, right=of 1, minimum size=1cm] (2) {2};
  \node[roundnode, right=of 2, minimum size=1cm] (3) {3};
  \node[roundnode, right=of 3, minimum size=1cm] (4) {4};
  \node[roundnode, right=of 4, minimum size=1cm] (5) {5};
  \node[roundnode, right=of 5, minimum size=1cm] (6) {6};

  \draw[solidarrow, thick] (1) to [bend left=40] node[below] {\bf 121}  (2);
  \draw[solidarrow, thick] (2) to [bend left=40] node[below] {\bf 137} (3);
  \draw[solidarrow, thick] (3) to [bend left=40] node[below] {\bf 75.4} (4);
  \draw[solidarrow] (4) to [bend left=40] node[below] {84.7} (5);
  \draw[solidarrow] (5) to [bend left=40] node[below] {78.5} (6);

  \draw[solidarrow] (2) to [bend left=40] node[below] {356} (5);
  \draw[solidarrow] (3) to [bend left=40] node[below] {238} (6);


  \draw[solidarrow] (1) to [bend right=40] node[below] {353} (3);
  \draw[solidarrow] (2) to [bend right=40] node[below] {225} (4);
  \draw[solidarrow] (3) to [bend right=40] node[below] {149} (5);
  \draw[solidarrow, thick] (4) to [bend right=40] node[below] {\bf 142} (6);


\end{tikzpicture}
}
\caption{Filtered graph structure of the relaxed problem for our numerical example, the shortest path is in bold}
\label{fig:shortestPathExample_filtered}
\end{figure}

\subsection{State Space Augmentation}

By applying a suitable algorithm (e.g. Dijkstra's algorithm) we find the shortest path from node 1 to node $T+1$ in the state space graph $G$.
As previously discussed, the optimal path can be immediately matched to a replenishment plan for problem $D$, and feasibility of this plan (i.e. absence of expected negative orders) can also be readily checked. If the plan is feasible w.r.t. problem $D$, then it is optimal. Conversely, if at any period $t$ we observe that $S_t < \tilde{I}_{t-1}$, this means an expected negative order is scheduled, and this infeasibility must be rectified to ensure the replenishment plan is feasible w.r.t. problem $D$.

In this subsection, we describe the high-level procedure for solving the non-stationary stochastic lot sizing problem using our SSA approach. After constructing the directed graph $G$ representing the state space, we solve a shortest path problem on this graph, and then iteratively fix any negative orders detected by modifying the graph. The pseudocode of the method is provided in Algorithm~\ref{alg:SSA-solve}.

\begin{algorithm}[ht]
\caption{State Space Augmentation}
\label{alg:SSA-solve}
\begin{algorithmic}[1]
\State \textbf{Input:} an instance of the Stochastic Lot Sizing problem
\State Create the state space graph: $G$ $\gets$ \texttt{create\_di\_graph()}
\State $p$ $\gets$ shortest path from node 1 to node $T+1$ in $G$
\While{the optimal path contains negative orders}
    \State Detect negative orders and fix them via procedure \texttt{fix\_negative\_orders(p,G)}, obtaining $G'$
    \If{at least one negative order has been detected and fixed}
        \State $p$ $\gets$ shortest path from node 1 to node $T+1$ in $G'$
        \State $G$ $\gets$ $G'$
    \EndIf
\EndWhile
\State Convert the final path to a binary representation and retrieve the order-up-to levels
\State \textbf{Output:} The optimal cost, binary order decisions, and order-up-to levels.
\end{algorithmic}
\end{algorithm}

The algorithm begins by generating a directed graph via procedure \texttt{create\_di\_graph()}. It then finds the shortest path from node 1 to node $T+1$ in $G$.
Next, the algorithm checks for expected negative orders along the identified path. If occurrences are detected, the method \texttt{fix\_negative\_orders}, discussed in the next section, is invoked to modify the graph by augmenting the state space. This fix-and-resolve procedure is performed iteratively until a path free of negative orders is obtained. Finally, the binary order decisions and the corresponding order-up-to levels are extracted and returned together with the optimal cost.

\subsection{Fixing Expected Negative Orders}

In our approach, a negative expected order arises when the expected closing inventory of one cycle exceeds the order-up-to level of the following cycle, which violates constraint \eqref{eq:DE-Q-2}. This means the two cycles cannot be optimised independently. 
Exploiting the convexity of the multiperiod newsvendor cost function, the infeasibility can be resolved by merging the cycles into a single longer cycle and optimising to obtain a single order-up-to level for the merged cycles. Our approach achieves this by (i) adding new nodes to the graph, (ii) reconnecting these nodes to repair the expected negative order transitions, and (iii) attaching the repaired segment to future nodes. These steps are explained below, and the detailed pseudocode is given in Algorithm \ref{alg:fix_neg_orders}.

First, the algorithm initialises variables and constructs a list of cycles from the input path (Lines 1--13 in Algorithm \ref{alg:fix_neg_orders}). It then detects the first occurrence of a negative order (Line 8) and defines the range of cycles affected (Lines 9--14). If no negative orders are found, the current graph is returned (Line 15).

\begin{algorithm}[p]
\caption{Fix Negative Orders}
\label{alg:fix_neg_orders}
\begin{algorithmic}[1]
\Procedure{fix\_negative\_orders}{path, $G$}
    \State \textbf{Initialize:}
        \Statex negative\_order\_fixed $\gets$ \textbf{False}
        \Statex pred $\gets$ \textbf{None}
        \Statex cycles $\gets$ empty list
    \For{each \textit{node} in path} \Comment{Construct the cycle list}
        \If{pred is \textbf{None}}
            \State pred $\gets$ \textit{node}
        \Else
            \State Retrieve arc from $G$ between pred and \textit{node}
            \State Append arc to cycles
            \State pred $\gets$ \textit{node}
        \EndIf
    \EndFor
    \For{$i = 0$ to $|\text{cycles}| - 2$} \Comment{Detect first negative order}
        \If{cycles[$i$].$\tilde{I}$ $>$ cycles[$i+1$].$S$}
            \State start $\gets i$
            \State \textbf{break}
        \EndIf
    \EndFor
    \If{start is defined} \Comment{Detect end of negative order cycle sequence}
        \If{there exists an index $j >$ start such that cycles[$j$].$\tilde{I}$ $<$ cycles[$j$].$S$}
            \State end $\gets j$
            \State negative\_order\_fixed $\gets$ \textbf{True}
        \Else
            \State end $\gets |\text{cycles}| - 1$
            \State negative\_order\_fixed $\gets$ \textbf{True}
        \EndIf
    \EndIf
    \If{negative\_order\_fixed is \textbf{False}} \Comment{No negative orders found}
        \State \Return (negative\_order\_fixed, $G$)
    \Else \Comment{Augment the state space}
        \State new\_nodes $\gets$ \Call{AddNewNodes}{$G$, start, end, cycles} \label{alg:line_newnodes_new}
        \State Remove the arc corresponding to the first negative order at index \texttt{start} from $G$
        \State mpnbs $\gets$ list of multiperiod newsvendor models for cycles from \texttt{start} to \texttt{end} \label{alg:line_mpnbs_new}
        \State Define function $f$, which aggregates the costs of mpnbs for a given order-up-to level S \label{alg:def_f}
        \State Compute the optimal order-up-to level S by minimizing $f$ \label{alg:min_f}
        \State \Call{ConnectSequentially}{$G$, new\_nodes, cycles, start, end, $S$} \label{alg:line_connect_new}
        \State \Call{ConnectToSubsequentNodes}{$G$, new\_nodes[0], cycles, end} \label{alg:line_connect_subseq}
        \State \Return (negative\_order\_fixed, updated $G$)
    \EndIf
\EndProcedure
\end{algorithmic}
\end{algorithm}

When a negative order is detected, the algorithm repairs it by first adding new nodes for the affected cycles via the \texttt{AddNewNodes} procedure (Line \ref{alg:line_newnodes_new} and Algorithm \ref{alg:helper_procedures}). Next, after removing the problematic arc, the multiperiod newsvendor models corresponding to the affected cycles are gathered (Line \ref{alg:line_mpnbs_new}). The algorithm then defines a new combined expected total cost function $f$ that aggregates the expected total costs of the multiperiod newsvendor models, by accounting for the expected reduction in available inventory caused by previous cycles' expected demand (Line \ref{alg:def_f}). More formally, consider two cycles $R(i,k-1)$ and $R(k,j-1)$, let $S_k$ be the order-up-to level of $R(k,j-1)$, and $\tilde{I}_{k-1}$ be the expected closing inventory level of $R(i,k-1)$. Recall that in model $D$, if $\delta_t=1$ then $S_t=\tilde{d}_t+\tilde{I}_t$; an expected negative order at period $k$ entails a violation of Eq. \ref{eq:DE-Q-2}, that is $\tilde{I}_{k-1}\geq S_k$. To resolve this violation, note that due to the convexity of replenishment cycle expected total cost functions, when negative order quantity scenarios arise, at optimality the expected closing inventory level of the first cycle must be equal to the opening inventory level of the second cycle \citep[see][pp. 237-238]{Rossi2007}. We therefore combine the expected total cost functions $C_{i,k}(\cdot)$ and $C_{k,j}(\cdot)$ associated with the multiperiod newsvendor models of $R(i,k-1)$ and $R(k,j-1)$ via the linking constraint $\tilde{I}_{k-1}=S_k$, and obtain a new expected total cost function $C^*_{i,j}(S)=C_{i,k}(S)+C_{k,j}(S-\tilde{d}_i-\ldots-\tilde{d}_{k-1})$. This function is minimised to compute an overall order-up-to level $S$ for the sequence of affected cycles (Line \ref{alg:min_f}). The new nodes are then connected in sequence using the \texttt{ConnectSequentially} procedure (Line \ref{alg:line_connect_new} and Algorithm \ref{alg:helper_procedures}), and finally, a bridging node is created (inheriting the label of the first repaired cycle) and is linked forward to all future nodes and backward to the node preceding the repaired segment to ensure continuity via \texttt{ConnectToSubsequentNodes} (Line \ref{alg:line_connect_subseq} and Algorithm \ref{alg:helper_procedures}). The updated graph is then returned.

\begin{algorithm}[p]
\caption{Sub-procedures}
\label{alg:helper_procedures}
\begin{algorithmic}[1]
\Procedure{AddNewNodes}{G, start, end, cycles}
    \State Initialize empty list new\_nodes
    \For{$i = $ start to end $- 1$}
        \State new\_node $\gets$ $|V|+1$, where $V$ is the set of nodes in $G$
        \State Add new\_node to $G$ with label equal to cycles[$i+1$].source
        \State Append new\_node to new\_nodes
    \EndFor
    \State \Return new\_nodes
\EndProcedure
\end{algorithmic}
\vspace{1em}
\begin{algorithmic}[1]
\Procedure{ConnectSequentially}{G, new\_nodes, cycles, start, end, S}
    \State current\_node $\gets$ cycles[start].source
    \For{$i = 0$ to $|\text{new\_nodes}| - 1$}
        \State Add an arc $a$ from current\_node to new\_nodes[$i$] in $G$ with:
            \Statex \quad \quad \quad $a.C$ = expected total cost computed for cycles[$i$] using order-up-to level $S$
            \Statex \quad \quad \quad $a.S$ = S
            \Statex \quad \quad \quad $a.\tilde{I}$ = S minus the expected demand of cycles[$i$]
        \State Update S $\gets$ S minus the expected demand in cycles[$i$]
        \State current\_node $\gets$ new\_nodes[$i$]
    \EndFor
    \State Add an arc $a$ from current\_node to cycles[end].dest in $G$ with:
         \Statex \quad \quad \quad $a.C$ = expected total cost computed for cycles[end] using order-up-to level S
         \Statex \quad \quad \quad $a.S$ = S
         \Statex \quad \quad \quad $a.\tilde{I}$ = S minus the expected demand of cycles[end]
\EndProcedure
\end{algorithmic}
\vspace{1em}
\begin{algorithmic}[1]
\Procedure{ConnectToSubsequentNodes}{G, first\_new\_node, cycles, start}
    \State Let $T$ be the total number of periods and let $V$ denote the set of nodes in $G$.
    \State new\_aux\_node $\gets |V|$
    \State Add node new\_aux\_node to $G$ with label equal to cycles[start+1].source.
    \State Let $q \gets$ label of cycles[start+1].source.
    \State Let $p \gets$ label of cycles[start+1].dest.
    \For{each arc $(u,v)$ in $G$ such that $u = q$  and $label(v) > p$}
         \State Duplicate arc $(u,v)$ by adding an arc from new\_aux\_node to $v$ with identical attributes.
    \EndFor
    \State $s\gets$ label of cycles[start].source.
    \State Construct multiperiod newsvendor model $m$ for cycle $R(s,q-1)$.
    \State Add an arc from cycles[start].source to new\_aux\_node in $G$
          \Statex \quad \quad \quad $a.C = C_{s,q}(S)$
          \Statex \quad \quad \quad $a.S = m.S$
          \Statex \quad \quad \quad $a.\tilde{I} = m.S$ minus expected demand of cycle $R(s,q-1)$.
\EndProcedure
\end{algorithmic}
\end{algorithm}

This augmentation procedure ensures that the final solution satisfies constraint \ref{eq:DE-Q-2}, which forces expected orders to be nonnegative; and thus, it ensures feasibility of replenishment plans generated w.r.t. model $D$.

\paragraph{Numerical example.}  After applying procedure \texttt{fix\_negative\_orders} to the shortest path previously obtained, we obtain the augmented graph $G'$ illustrated in Figure \ref{tab:G_1} and Figure \ref{fig:shortestPathExample_augmented}; note that new nodes 7 and 8 have label 3, to indicate that these new nodes have been created to duplicate existing node 3, this label is necessary to reconstruct the optimal replenishment plan from the shortest path. For $G'$, we obtain a shortest path (shown in bold in Figure \ref{tab:G_1}) with cost 486.5 and the replenishment plan illustrated in Figure \ref{fig:replenishment_plan_2}, for which no further negative orders were detected, and which is therefore optimal for problem $D$.

\begin{figure}[t]
\begin{minipage}[b]{0.45\columnwidth}
\begin{center}
\resizebox{.85\textwidth}{!}{
    \begin{tabular}{ccccccccc}
        \toprule
        & 1 & 2 & 3 & 4 & 5 & 6 & $7_3$ & $8_3$\\
        \midrule
1 & -- & 121 & 353 & -- & -- & -- & -- & -- \\
2 & -- & -- & -- & 225 & 356 & -- & 143 & 137 \\
3 & -- & -- & -- & 75.4 & 149 & 238 & -- & -- \\
4 & -- & -- & -- & -- & 84.7 & 142 & -- & -- \\
5 & -- & -- & -- & -- & -- & 78.5 & -- & -- \\
6 & -- & -- & -- & -- & -- & -- & -- & -- \\
$7_3$ & -- & -- & -- & 83.30 & -- & -- & -- & -- \\
$8_3$ & -- & -- & -- & -- & 150 & 238 & -- & -- \\
        \bottomrule
    \end{tabular}
}
\caption{Augmented connection matrix}
\label{tab:G_1}
\end{center}
\end{minipage}
\hfill
\begin{minipage}[b]{0.45\columnwidth}
\begin{center}
\resizebox{.75\textwidth}{!}{
\begin{tikzpicture}
    \begin{axis}[
        xlabel={Period},
        ylabel={inventory level},
        xmin=0, xmax=5,
        ymin=-50, ymax=200,
        xtick={0,1,...,5},
        xticklabels=\empty,
        extra x ticks={0.5,1.5,2.5,3.5,4.5,5.5},
        extra x tick labels={1,2,3,4,5},
        extra x tick style={major tick length=0pt},
        ytick={0,50,...,200},
        legend pos=north east,
        ymajorgrids=true,
        grid style=dashed,
    ]
    \addplot[
        color=blue,
        mark=x,
        ]
        coordinates {
        (0,149)(1,149-100)(1,187)(2,187-125)(2,83)(3,83-25)(4,83-25-40)(4,45)(5,45-30)
        };
    \addplot[
        color=red,
        only marks,
        mark=*,
        ]
        coordinates {
        (0,149) (1,187) (2, 83) (4, 45)
        };
    \end{axis}
\end{tikzpicture}
}
\caption{Replenishment plan after first iteration}
\label{fig:replenishment_plan_2}
\end{center}
\end{minipage}
\end{figure}

\begin{figure}[b]
\centering
\resizebox{0.75\textwidth}{!}{%
\begin{tikzpicture}[
  roundnode/.style={circle, draw, minimum size=0.6cm, inner sep=1pt},
  dashedarrow/.style={dashed, ->, -Latex},
  solidarrow/.style={->, -Latex},
  node distance=1.2cm and 1.5cm
]
  \node[roundnode, minimum size=1cm] (1) {1};
  \node[roundnode, right=of 1, minimum size=1cm] (2) {2};
  \node[roundnode, right=of 2, minimum size=1cm] (3) {3};
  \node[roundnode, right=of 3, minimum size=1cm] (4) {4};
  \node[roundnode, below=of 3, minimum size=1cm] (7) {$7_3$};
  \node[roundnode, below=of 7, minimum size=1cm] (8) {$8_3$};
  \node[roundnode, right=of 4, minimum size=1cm] (5) {5};
  \node[roundnode, right=of 5, minimum size=1cm] (6) {6};

  \draw[solidarrow, thick] (1) to [bend left=40] node[below] {\bf 121}  (2);
  \draw[solidarrow, thick, red] (2) to [bend left=40, red] node[below] {\bf 137} (3);
  \draw[solidarrow] (3) to [bend left=40] node[below] {75.4} (4);
  \draw[solidarrow] (4) to [bend left=40] node[below] {84.7} (5);
  \draw[solidarrow, thick] (5) to [bend left=40] node[below] {\bf 78.5} (6);

  \draw[solidarrow] (2) to [bend left=40] node[below] {356} (5);
  \draw[solidarrow] (3) to [bend left=40] node[below] {238} (6);


  \draw[solidarrow] (1) to [bend right=40] node[below] {353} (3);
  \draw[solidarrow] (2) to [bend right=40] node[below] {225} (4);
  \draw[solidarrow] (3) to [bend right=40] node[below] {149} (5);
  \draw[solidarrow] (4) to [bend right=40] node[below] {142} (6);


\draw[solidarrow] (2) to [bend right=40] node[below] {143} (7);
\draw[solidarrow] (7) to [bend right=40] node[below] {83.3} (4);

\draw[solidarrow, thick] (2) to [bend right=40] node[below] {\bf 137} (8);
\draw[solidarrow, thick] (8) to [bend right=40] node[below] {\bf 150} (5);
\draw[solidarrow] (8) to [bend right=40] node[below] {238} (6);

\end{tikzpicture}
}
\caption{Augmented graph structure of the relaxed problem for our numerical example, the shortest path is in bold, in red we identify arcs that have been removed by  \texttt{fix\_negative\_orders}}
\label{fig:shortestPathExample_augmented}
\end{figure}

\newpage

\section{Computational study}\label{sec:computation}


This section presents a computational study that investigates the performance of our newly proposed approach. All computational experiments are carried out on an Intel(R) Core(TM) i7-8650U CPU 16GB RAM; our approach is coded in Python 3.11; MILP problems are solved using Gurobi Optimiser v10.0.2.

As benchmark, we consider the cut generation (CG) method introduced by \citet{tunc2018extended}, which represents the fastest known approach for modelling and solving non-stationary $(R,S)$ policies under the penalty cost scheme. Our experiments aim to demonstrate that our method outperforms this benchmark in terms of computational efficiency. 

It should be emphasised that both approaches solve the deterministic non-linear model $D$ presented in Section \ref{sec:DE} to optimality. \citet{tunc2018extended} relies on tight linear cuts, while our approach uses numerical optimisation to optimise the multiperiod newsvendor models in the state-space graph; therefore, in our computational study we will not discuss optimality gaps. 

Our test bed is based on Set-A of \citet{tunc2018extended}, which assumes normally distributed demand with a fixed coefficient of variation, and includes planning horizons $T \in \{40, 100, 250\}$, fixed ordering cost $K \in \{225, 900, 2500\}$, penalty parameters $b \in \{2, 5, 10\}$, and coefficients of variation $\rho \in \{0.1, 0.2, 0.3\}$, where $\rho=\sigma_t/\mu_t$. 
For the erratic demand pattern, mean demands per period are uniformly distributed as $\tilde{d}_t \sim \mathcal{U}(2,100)$, and for the lumpy demand pattern, demand is $\tilde{d}_t \sim \mathcal{U}(2,420)$ with probability 0.2 and $\tilde{d}_t \sim \mathcal{U}(2,20)$ with probability 0.8. 
We note that the lower bound of the uniform distribution in \citeauthor{tunc2018extended}'s test set is originally 0; however, in our study we raise this to 2 to avoid numerical instability and division-by-zero issues that may arise when the mean of the distribution is close to zero. 
We generated 10 independent problem instances for each demand pattern, and thus obtained 1620 instances in total. 

Since demand is assumed to be normally distributed, in this computational study the computation of the first order loss function and of its complementary function in $C_{(i,j)}(S)$ is carried out by using  relevant closed form expressions \cite[see][Section 3]{rossi2014piecewise}; similarly, these closed form expressions have been used in our re-implementation of \cite{tunc2014reformulation}.

{\footnotesize{
\begin{table}[H]\renewcommand{\arraystretch}{0.8}
\caption{Pivot table illustrating runtimes of SSA and CG}\label{tab:compare}
\centering
\begin{tabular}{@{}
  >{\raggedleft\arraybackslash}p{1.1cm}
  >{\raggedleft\arraybackslash}p{1.1cm}|
  >{\centering\arraybackslash}p{2.5cm}
  >{\centering\arraybackslash}p{2cm}|
  >{\centering\arraybackslash}p{2.5cm}
  >{\centering\arraybackslash}p{2.5cm}
@{}}\toprule
            &values       &SSA runtime average &CG runtime average     &SSA runtime median &CG runtime median\\
\midrule
pattern     &Erratic        &59.1      &2321      &30.1        &90.0\\
            &Lumpy          &111     &1857      &60.4        &86.9\\
\midrule
$T$         &40     &13.6      &2.15      &12.1       &1.92\\
            &100    &49.2      &101      &47.1     &90.3\\
            &250    &192        &6163    &155     &6171\\
\midrule
$K$         &225    &32.3      &1401      &18.9        &61.9\\
            &900    &75.3      &2064      &42.8        &89.7\\
            &2500   &147     &2800      &72.0        &106\\
\midrule
$b$         &2      &90.2      &1396      &40.1        &79.9\\
            &5      &83.4      &2065      &36.1        &91.4\\
            &10     &81.1      &2804      &31.8        &98.2\\
\midrule
$\rho$      &0.1    &86.4      &1393      &36.8        &75.2\\
            &0.2    &84.3      &2065      &36.6        &93.0\\
            &0.3    &83.9      &2837      &37.2        &112\\
\midrule
            &Total      &84.9      &2088      &36.8        &90.3\\
\bottomrule
\end{tabular}
\end{table}
}}

Finally, because \citeauthor{tunc2018extended}'s method becomes prohibitively slow for long planning horizons, we do not attempt to run it on all 540 250-period instances in our test bed. Instead, we manually sampled a small but representative subset of $N = 12$ problem instances: for each pattern, two independent instances at each of the three $(K,b,\rho)$ parameter settings, which are drawn from the same $T=250$ problem pool used to evaluate our SSA approach, and \citeauthor{tunc2018extended}'s results at $T=250$ are reported as the mean runtime over these sampled runs.
In the reported results, the mean CG runtime corresponds to the mean of the means across the three $T$ settings. Additionally, we report the median CG runtime, which is obtained by bootstrapping the 12 sampled runtimes to form 540 instances for $T=250$; therefore, three sets of $T$ are balanced, from which the median is calculated.




\begin{figure}[H]
    \centering
    \includegraphics[width=0.8\linewidth]{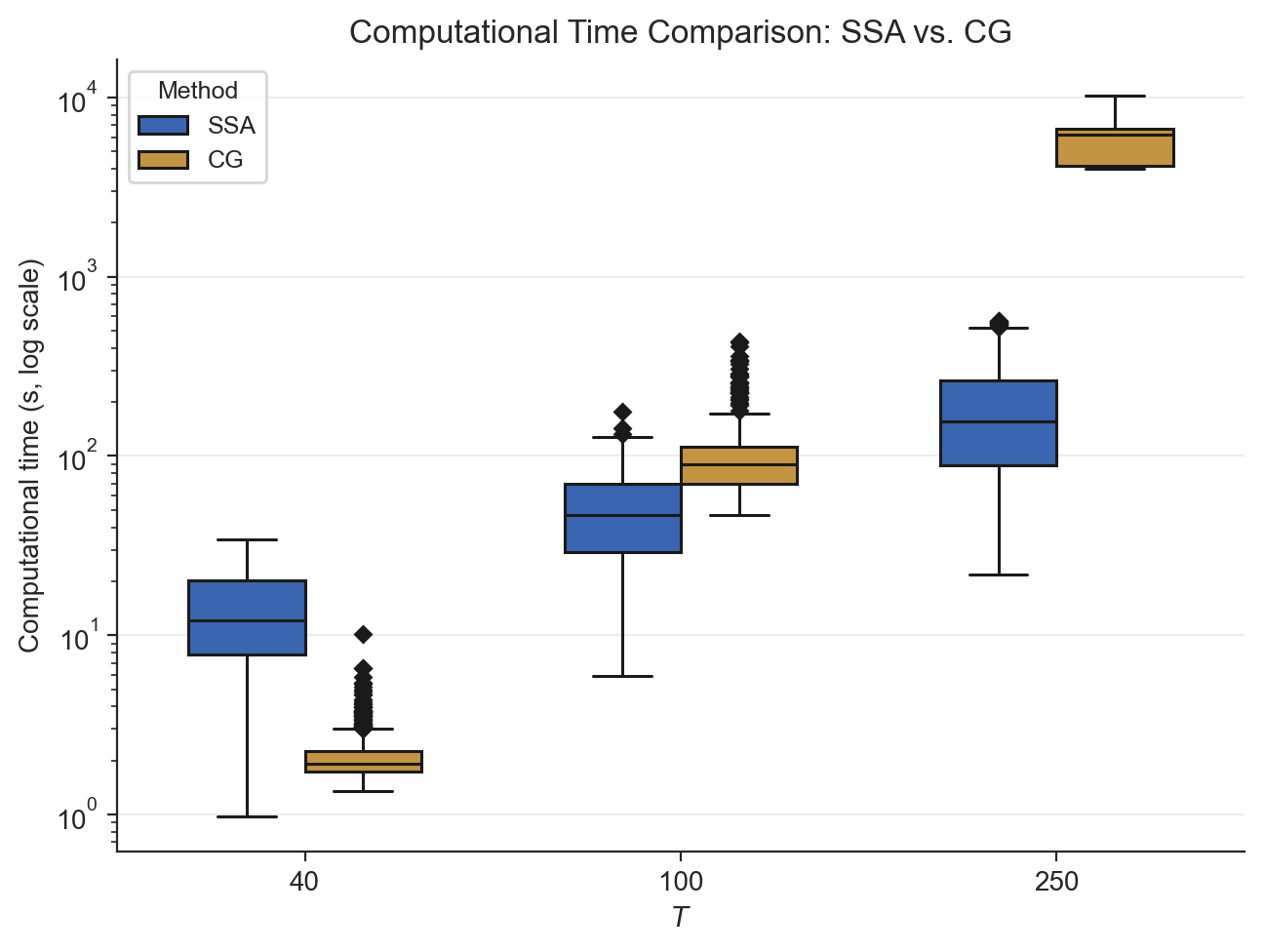}
    \caption{Computational Time Comparison: SSA vs. CG.}
    \label{fig:saa_cg_compare_bootstrapping}
\end{figure}

Figure~\ref{fig:saa_cg_compare_bootstrapping} shows the distribution of the computational time of SSA and CG over the test instances at each planning horizon $T$. The vertical axis is logarithmic. At $T=40$ the exact CG benchmark is faster because the instances are small, but SSA overtakes it at $T=100$ and is far faster at $T=250$, where the run time of CG grows steeply.
At every horizon, the slowest SSA runs grow only gradually with $T$ and never reach the extreme times that CG exhibits for long horizons, which shows that SSA stays robust across the instance space.

Table~\ref{tab:compare} summarises the average and median computational times of SSA and \citeauthor{tunc2018extended}'s CG method when pivoting on each parameter. As shown, SSA’s average runtimes differ substantially by demand pattern (59.1s vs 111s), with corresponding medians of 30.1s and 60.4s; in contrast, CG’s averages are much higher (2321s and 1857s) with its medians also higher (90.0s and 86.8s). The fact that the medians for CG are much lower than the averages implies that most runs remain tractable despite a few extremely slow instances for long planning horizons.
Both methods slow as the horizon $T$ increases, but CG’s distribution is far more skewed: at $T=250$ its average time soars to 6163s (median 6171s), whereas SSA’s mean and median are 192s and 155s. As the order-up-to level $K$ grows, SSA runtime rises from 32.3s to 147s , while CG increases from 1401s to 2800s.
When the penalty cost $b$ increases, SSA’s mean decreases slightly (from 90.2s to 81.1s) with medians falling from 40.1s to 31.8s, whereas CG's runtime increases on average from 1396s to 2804s and its median rises from 79.9s to 98.2s. 
A similar pattern emerges for $\rho$, where SSA remains stable around 84s on average, but CG’s average runtime more than doubles and its median grows from 75.2s to 112s.
Overall, SSA solves all instances in 84.9s on average (median 36.8s), compared to 2088s for CG (median 90.3s).  
Our study underscores SSA’s ability to solve long horizon problems efficiently, whereas CG becomes increasingly intractable as $T$ grows.

For completeness, Appendix B reports the average SSA augmented graph sizes and CG cut counts over pivots on $K$, $b$, and $\rho$ for each horizon $T$. While per‐arc and per‐cut computational costs are not directly comparable, Table \ref{tab:iterations} shows the dramatic impact of the state space filtering discussed in Section \ref{sec:state_space_filtering} and that the augmented state space graph remains very compact for all instances. In contrast, the number of cuts generated by CG displays a substantial growth with $T$.

\section{Conclusion}\label{sec:conclusion}

In this paper, we developed a relaxation-and-augmentation algorithm 
to compute near-optimal $(R,S)$ policy parameters for the single‐item single-stocking location non-stationary stochastic lot‐sizing problem under penalty costs. Our approach transforms the relaxed stochastic lot-sizing problem into a shortest-path problem on a filtered directed acyclic graph, and then repairs infeasible negative-order states via a repetitive state space augmentation procedure. By exploiting dominance relationships and convexity of the multi-period newsvendor subproblems, our method prunes irrelevant arcs from the state space graph and dramatically reduces the search effort.

Comprehensive experiments on 1620 benchmark instances demonstrate that our approach outperforms the state-of-the-art CG approach of \citet{tunc2018extended}: on medium-scale problems our method completes search in 50s on average, while CG requires 100s on average; 
for very long horizons ($T=250$), our method completes search in 192s on average, while CG requires in excess of 6000s per instance.
The performance of our approach is consistent across erratic and lumpy demand, varying ordering and penalty costs, and coefficients of variation. These results demonstrate the scalability of our approach. Future research can focus on deriving worst-case complexity bounds and extending the approach to $(s,Q)$ and $(s,S)$ policies.  

\clearpage

\appendix


\renewcommand\thefigure{\Alph{section}\arabic{figure}}
\renewcommand\theequation{\Alph{section}\arabic{equation}}
\setcounter{figure}{0}
\setcounter{table}{0}
\setcounter{equation}{0}
\renewcommand\thetable{\Alph{section}\arabic{table}}

\renewcommand\thefigure{\Alph{section}\arabic{figure}}
\setcounter{figure}{0}
\setcounter{equation}{0}
\section{Notations}\label{sec:App.notations}

{\small{
\begin{longtable}{rp{12cm}}
\caption{Notation summary}\label{tab:App.notations}\\
\toprule
Notation & Description \\
\midrule
\endfirsthead
\multicolumn{2}{c}{\tablename~\thetable{} (continued)}\\
\toprule
Notation & Description \\
\midrule
\endhead
\midrule \multicolumn{2}{r}{Continued on next page}\\
\endfoot
\bottomrule
\endlastfoot

$T$ & planning horizon length (number of periods).\\
$d_t$ & stochastic demand in period $t$.\\
$\tilde d_t$ & expected demand in period $t$, i.e.\ $\mathrm{E}[d_t]$.\\
$d_{i,k}$ & cumulative demand from period $i$ to $k$, $d_{i,k}=\sum_{r=i}^k d_r$.\\
$\delta_t$ & binary indicator, equals 1 if an order is placed at the start of period $t$, 0 otherwise.\\
$S_t$ & order‐up‐to level for the cycle beginning at period $t$.\\
$Q_t$ & order quantity in period $t$, given by $Q_t=\delta_t\max\{S_t-I_{t-1},0\}$.\\
$I_t$ & inventory level at the end of period $t$ (negative values denote backorders).\\
$L_t(y)$ & immediate expected holding and penalty cost when post‐replenishment inventory is $y$:
$L_t(y)=h\,\mathrm{E}[\max(y-d_t,0)] + b\,\mathrm{E}[\max(d_t-y,0)]$.\\
$\tilde H_t(y)$ & expected on‐hand inventory at end of period $t$: $\mathrm{E}[\max(y-d_t,0)]$.\\
$\tilde B_t(y)$ & expected backorders at end of period $t$: $\mathrm{E}[\max(d_t-y,0)]$.\\
$\tilde I_t$ & expected closing inventory at end of period $t$, $\tilde I_t=\tilde H_t(S_t)-\tilde B_t(S_t)$.\\
$C_1(I_0)$ & minimum expected total cost from period 1 to $T$ given initial inventory $I_0$.\\
$C_{i,j}(S)$ & expected cost of replenishment cycle $R(i,j-1)$ with order‐up‐to level $S$, where $C_{i,j}(S)=K + \sum_{k=i}^{j-1}\Bigl[h\,\mathrm{E}[\max(S-d_{i,k},0)]+b\,\mathrm{E}[\max(d_{i,k}-S,0)]\Bigr]$.\\
$G=(V,E)$ & directed acyclic graph whose nodes $V=\{1,\dots,T+1\}$ represent review epochs (and terminal), and whose arcs $E$ represent replenishment cycles.\\
$V$ & set of nodes (period start and dummy terminal).\\
$E$ & set of arcs $(i,j)$, each denoting cycle $R(i,j-1)$ with cost $C_{i,j}(S_i)$.\\
$(i,j)$ & arc from node $i$ to node $j$, corresponding to cycle $R(i,j-1)$.\\
$p$ & shortest path from node 1 to node $T+1$ in $G$, identifying the optimal set of cycles under relaxation.\\

\end{longtable}
}}

\renewcommand\thefigure{\Alph{section}\arabic{figure}}
\renewcommand\theequation{\Alph{section}\arabic{equation}}
\setcounter{figure}{0}
\setcounter{table}{0}
\setcounter{equation}{0}
\renewcommand\thetable{\Alph{section}\arabic{table}}

\renewcommand\thefigure{\Alph{section}\arabic{figure}}
\setcounter{figure}{0}
\setcounter{equation}{0}
\section{Averaged Computational Metrics by Pivot Parameters}\label{sec:App.graphCuts}

{\scriptsize{
\begin{table}[H]\renewcommand{\arraystretch}{0.8}
\caption{Computational statistics for SSA (average augmented graph size) and CG (average number of cuts) pivoting parameters under $T$ setup.}\label{tab:iterations}
\centering
\begin{tabular}{@{}
  >{\raggedleft\arraybackslash}p{0.6cm}
  >{\raggedleft\arraybackslash}p{0.8cm}|
  >{\centering\arraybackslash}p{2.6cm}
  >{\centering\arraybackslash}p{1cm}|
    >{\centering\arraybackslash}p{2.6cm}
  >{\centering\arraybackslash}p{1cm}|
    >{\centering\arraybackslash}p{2.6cm}
  >{\centering\arraybackslash}p{1cm}
@{}}\toprule
&&$\qquad\qquad T=40$&&$\qquad\quad T=100$&&$\qquad\quad T=250$&\\
\midrule
        &values       &{\small{SSA: augmented graph size}}     &{\small{CG: cuts}}        &{\small{SSA: augmented graph size}}          &{\small{CG: cuts}}       &{\small{SSA: augmented graph size}}      &{\small{CG: cuts}}\\
\midrule
$K$         &225    &$\langle42,291\rangle$     &370    &$\langle104,857\rangle$    &1152   &$\langle258,2447\rangle$   &1024\\
            &900    &$\langle41,421\rangle$     &439    &$\langle101,1318\rangle$   &1344   &$\langle252,3835\rangle$   &3667\\
            &2500   &$\langle41,533\rangle$     &397    &$\langle101,1850\rangle$   &1469   &$\langle251,5497\rangle$   &13365\\
\midrule
$b$         &2      &$\langle41,446\rangle$     &262    &$\langle101,1434\rangle$   &885    &$\langle252,4221\rangle$   &1024\\
            &5      &$\langle41,417\rangle$     &394    &$\langle102,1318\rangle$   &1310   &$\langle254,3847\rangle$   &3667\\
            &10     &$\langle42,403\rangle$     &550    &$\langle103,1273\rangle$   &1813   &$\langle255,3711\rangle$   &13365\\
\midrule
$\rho$      &0.1    &$\langle41,422\rangle$     &153    &$\langle101,1340\rangle$   &453    &$\langle251,3913\rangle$   &1024\\
            &0.2    &$\langle41,422\rangle$     &372    &$\langle102,1340\rangle$   &1218	&$\langle253,3917\rangle$   &3667\\
            &0.3    &$\langle42,434\rangle$     &882    &$\langle102,1345\rangle$   &2337   &$\langle257,3949\rangle$   &13365\\
\midrule
    &total          &$\langle41,422\rangle$     &402    &$\langle102,1342\rangle$   &1336   &$\langle254,3926\rangle$   &6019\\
\bottomrule
\end{tabular}
\end{table}
}}

\clearpage
\bibliography{referenceReMILP}

\end{document}